\documentclass[12pt,a4paper]{article}
\usepackage[utf8]{inputenc}
\usepackage[T1]{fontenc}
\usepackage[english]{babel}
\usepackage{amsmath}
\usepackage{amsfonts}
\usepackage{amssymb}
\usepackage{graphicx}

\newtheorem{theorem}{Theorem}[section]
\newtheorem{corollary}{Corollary}[section]
\newtheorem{lemma}{Lemma}[section]

\begin{document}
	\begin{center}
		{\Large \textbf{A second order difference scheme for time fractional diffusion equation with generalized memory kernel}} 
		\vskip 5mm	{\normalsize Aslanbek Khibiev\footnote{Corresponding author}, Anatoly Alikhanov$^{2}$, Chengming Huang$^{3}$} 		
		\vskip 5mm
		
	\end{center}
		\begin{scriptsize}
			$^{1}$ \quad Institute of Applied Mathematics and Automation KBSC RAS,
			Nalchik 360000, Russia; akkhibiev@gmail.com \\
			$^{2}$ \quad North-Caucasus Center of Mathematical Research, North-Caucasus Federal University,   Stavropol 355017,  Russia; aaalikhanov@gmail.com\\
			$^{3}$ \quad School of Mathematics and Statistics, Huazhong University of Science and Technology, Wuhan 430074, China; chengming\_huang@hotmail.com
		\end{scriptsize}

	\abstract{In the current work we build a difference analog of the  Caputo
		fractional derivative with generalized memory kernel ($_\lambda$L2-1$_\sigma$ formula).
		The fundamental features of this difference operator are studied
		and on its ground some difference schemes generating approximations
		of the second order in time for the generalized time-fractional diffusion equation with variable coefficients are worked out. We have proved stability and convergence of the given schemes in the grid $L_2$ - norm with the rate equal to the order of the approximation error. The achieved results are supported by the numerical
		computations performed for some test problems.}
		
		
		\section{Introduction}
		\sloppy{ Differential equations with fractional order derivatives
			represent a powerful mathematical tool for exact and realistic
			description of physical and chemical processes for which it is
			needed to take into consideration the background (memory) of the
			process \cite{OldSpan,Podl,Hilfer,KilbSrivTruj}. The
			patterns which take memory into consideration in such equations are
			the memory functions that are the kernels of integrals defining the
			operators of fractional integro-differentiation. For fractional integro-differentiation
			operators, the memory functions are namely power
			functions. The exponent of the power function of memory defines
			the order of the derivative and is connected with the fractal dimension
			of the environment in which the described process takes place.
			For more accurate description of the process in heterogeneous porous media,
			differential equations with fractional derivatives of distributed
			order are often used too. Processes of the memory can be described with the help of the
			memory function which has more complex structure than the power
			function.}
		
		In the rectangle $\bar Q_T=\{0\leq x\leq 1, 0\leq t\leq T\}$ we consider
		the Dirichlet boundary value problem for time fractional diffusion
		equation with generalized memory kernel and variable coefficients
		\begin{equation}\label{ur1}
			\partial_{0t}^{\alpha,\lambda(t)}u=\mathcal{L}u+f(x,t)
			,\quad 0<x<1,\quad 0<t\leq T,
		\end{equation}
		\begin{equation}
			u(0,t)=0,\quad
			u(1,t)=0,\quad 0\leq t\leq T,\quad  u(x,0)=u_0(x),\quad 0\leq x\leq 1, \label{ur3}
		\end{equation}
		where
		$$
		\mathcal{L}u=\frac{\partial }{\partial x}\left(k(x,t)\frac{\partial
			u}{\partial x}\right)-q(x,t)u,
		$$
		$$
		\partial_{0t}^{\alpha,\lambda(t)}u(x,t)=\frac{1}{\Gamma(1-\alpha)}\int\limits_{0}^{t}{\frac{\lambda(t-\eta)}{(t-\eta)^{\alpha}}}\frac{\partial
			u}{\partial \eta}(x,\eta)d\eta
		$$
		is the generalized Caputo fractional
		derivative of order $\alpha$, $0<\alpha<1$ with  weighting function
		$\lambda(t)\in \mathcal{C}^2[0,T]$, where $\lambda(t)>0$,
		$\lambda'(t)\leq0$ for all $t\in [0,T]$; $0<c_1\leq k(x,t)\leq c_2$,
		$q(x,t)\geq0$ for all $(x,t)\in \bar Q_T$.

		Diffusion and Fokker-Planck-Smoluchowski equations which have a generalized
		memory kernel were investigated in \cite{Fokk_Plan_Smol}. In this work
		it is demonstrated that the memory kernel appearing in the generalized
		diffusion equation has diverse potential forms which can describe a
		broad range of experimental phenomena.

		With the help of the energy inequality method, a priori estimates for the
		solution of both differential and difference problems of the
		Dirichlet and Robin boundary value problems for the fractional,
		variable and distributed order diffusion equation with Caputo
		fractional derivative were derived in \cite{Alikh:10, Alikh:12,
			AlikhanovJCP, Alikh_15, Alikh_16, Alikh_17_gen, Alikh_17, Khibiev}. A priori estimates for
		the difference problems analyzed in \cite{ShkhTau:06} by means of the
		maximum principle imply the stability and convergence of these
		difference schemes.
		
		{In this work, to construct difference schemes with the order
			of accuracy $ O(\tau^{2}) $ in time we have to demand the
			existence of a sufficiently smooth solution of the original problem.
			It brings on a significant narrowing of the input data class of the
			problem for which we apply the proposed method. As it is well known (see
			for example \cite{SakaYama,Luch}), in the case of smooth input data
			for a time-fractional diffusion equation, the solutions are not necessarily
			smooth in a closed domain, because the derivatives of the function
			$u(x,t)$ with respect to $t$ might possess a singularity at $t = 0$. In such cases, 
			if possible, we present the solution as the sum of
			two functions: one of which is known but not smooth, whereas the other
			is smooth but not known, as it is illustrated in work \cite{Alikh_17_1}.}
		
		{ In work \cite{Stynes}, we consider a reaction-diffusion problem with a Caputo
			time derivative of the order $\alpha\in(0,1)$. It is shown that
			the solution of such a problem has in general a weak
			singularity near the initial time $t = 0$, and we derive sharp pointwise
			bounds on certain derivatives of this solution. We have given a new
			analysis of a standard finite difference method for the problem,
			taking into account this initial singularity. }
		
		{In \cite{Lazarov_1}, we study an analysis of the L1 scheme for the subdiffusion equation with
			nonsmooth data. In \cite{Lazarov_2}}, error estimates
		for approximations of distributed order time fractional diffusion
		equation with nonsmooth data were investigated.
		
		In the current paper, a difference
		analog of the  Caputo fractional derivative with generalized memory
		kernel ($_\lambda$L2-1$_\sigma$ formula) is built up. The essential features of this difference
		operator are investigated and on its ground some difference schemes
		generating approximations of the second and fourth order in space
		and the second order in time for the generalized
		time-fractional diffusion equation with variable coefficients are
		studied. Stability of the suggested schemes as well as their
		convergence in the grid $L_2$ - norm with the rate equal to the
		order of the approximation error are proven. The achieved results
		are supported by the numerical computations performed for some
		test problems.

		\section{Stability and convergence of the family of difference schemes}
		
		In this section, we consider some families of difference schemes in a general form,
		set on a non-uniform time grid. A criterion of the
		stability of the difference schemes in the grid $L_2$ - norm is
		worked out. The convergence of solutions of the difference schemes to
		the solution of the corresponding differential problem with the rate
		equal to the order of the approximation error is proven.

		In the rectangle  $\overline Q_T=\{(x,t): 0\leq x\leq l,\, 0\leq
		t\leq T\}$ we assign the grid $\overline
		\omega_{h\tau}=\overline\omega_{h}\times\overline\omega_{\tau}$,
		where $\overline\omega_{h}=\{x_i=ih, \, i=0, 1, \ldots, N;\,
		hN=l\}$, $\overline\omega_{\tau}=\{t_j: \, 0=t_0<t_1<t_2<\ldots
		<t_{M-1}<t_{M}=T\}$.
		
		The family of difference schemes, approximating problem
		\eqref{ur1}--\eqref{ur3} on the grid $\overline \omega_{h\tau}$, mainly has
		the form
		
		\begin{equation}\label{ur03}
			{_g}\Delta_{0t_{j+1}}^{\alpha}y_i=\Lambda y^{(\sigma_{j+1})}_i
			+\varphi_i^{j+1}, \quad i=1,2,\ldots,N-1,\quad j=0,1,\ldots,M-1,
		\end{equation}
		\begin{equation}
			y(0,t)=0,\quad y(l,t)=0,\quad t\in \overline \omega_{\tau}, \quad
			y(x,0)=u_0(x),\quad  x\in \overline \omega_{h},\label{ur03.1}
		\end{equation}
		where
		\begin{equation}
			{_g}\Delta_{0t_{j+1}}^{\alpha}y_i=\sum\limits_{s=0}^{j}\left(y_i^{s+1}-y_i^s\right)g_{s}^{j+1},\quad
			g_{s}^{j+1}>0, \label{ur03.2}
		\end{equation}
		is a difference analog of the generalized Caputo fractional
		derivative of the order $\alpha$ with weighting function
		$\lambda(t)$ ($0<\alpha<1$, $\lambda(t)>0$, $\lambda'(t)\leq0$),
		$\Lambda$ is a difference operator which approximates the continuous
		operator $\mathcal{L}$, such that the operator $-\Lambda$ preserves
		its positive definiteness:
		$$
		(-\Lambda y,y)\geq \varkappa\|y\|_0^2, \quad
		(y,v)=\sum_{i=1}^{N-1}y_iv_ih, \quad \|y\|_0^2=(y,y), \quad
		\varkappa>0,
		$$
		$y^{j+\sigma_{j+1}}=\sigma_{j+1}y^{j+1}+(1-\sigma_{j+1})y^{j}$,
		$0\leq\sigma_{j+1}\leq1$, at $j=0,1,\ldots,M-1$.

		\begin{lemma}\label{lem_JCP}  \cite{AlikhanovJCP} If
			$g_{j}^{j+1}>g_{j-1}^{j+1}>\ldots>g_{0}^{j+1}>0$, $j=0,1,\ldots,M-1$
			then for any function $v(t)$ defined on the grid $\overline
			\omega_{\tau}$ the following inequalities hold true
			\begin{equation}\label{ur04}
				v^{j+1}{_g}\Delta_{0t_{j+1}}^{\alpha}v\geq
				\frac{1}{2}{_g}\Delta_{0t_{j+1}}^{\alpha}(v^2)+\frac{1}{2g^{j+1}_j}\left({_g}\Delta_{0t_{j+1}}^{\alpha}v\right)^2,
			\end{equation}
			\begin{equation}\label{ur05}
				v^{j}{_g}\Delta_{0t_{j+1}}^{\alpha}v\geq
				\frac{1}{2}{_g}\Delta_{0t_{j+1}}^{\alpha}(v^2)-\frac{1}{2\left(g^{j+1}_j-g^{j+1}_{j-1}\right)}\left({_g}\Delta_{0t_{j+1}}^{\alpha}v\right)^2,
			\end{equation}
			where $g^{1}_{-1}=0$.
		\end{lemma}

		\begin{corollary}\label{cor_JCP} 
			\cite{AlikhanovJCP} If
			$g_{j}^{j+1}>g_{j-1}^{j+1}>\ldots>g_{0}^{j+1}>0$ and
			$\frac{g_{j}^{j+1}}{2g_{j}^{j+1}-g_{j-1}^{j+1}}\leq\sigma_{j+1}\leq1$,
			where $j=0,1,\ldots,M-1$, $g_{-1}^{1}=0$, then for any function
			$v(t)$ defined on the grid $\overline\omega_{\tau}$ we have the
			inequality
			\begin{equation}\label{ur07}
				\left(\sigma_{j+1} v^{j+1}+(1-\sigma_{j+1})v^{j}\right){_g}\Delta_{0t_{j+1}}^\alpha v \geq \frac{1}{2}{_g}\Delta_{0t_{j+1}}^\alpha
				(v^2).
			\end{equation}
		\end{corollary}
		
		\begin{theorem}\label{theor_JCP_1} \cite{AlikhanovJCP} If	$$
			g_{j}^{j+1}>g_{j-1}^{j+1}>\ldots>g_{0}^{j+1}\geq c_2>0, \quad
			\frac{g_{j}^{j+1}}{2g_{j}^{j+1}-g_{j-1}^{j+1}}\leq\sigma_{j+1}\leq1,
			$$
			where  $j=0,1,\ldots,M-1$, $g_{-1}^{1}=0$, then the difference
			scheme \eqref{ur03}--\eqref{ur03.1} is unconditionally stable and
			its solution satisfies the following a priori estimate:
			\begin{equation}\label{ur08}
				\|y^{j+1}\|_0^2\leq\|y^0\|_0^2+\frac{1}{2\varkappa c_2}\max\limits_{0\leq j\leq
					M}\|\varphi^{j}\|_0^2,
			\end{equation}
		\end{theorem}

		A priori estimate  \eqref{ur08} implies the stability of difference
		scheme \eqref{ur03}--\eqref{ur03.1}.
		
		\begin{theorem}\label{theor_JCP_2} {\cite{AlikhanovJCP} If the conditions of
				Theorem \eqref{theor_JCP_1} are fulfilled and difference scheme
				\eqref{ur03}--\eqref{ur03.1} has the approximation order
				$\mathcal{O}(N^{-r_1}+M^{-r_2})$, where $r_1$ and $r_2$ are some
				known positive numbers, then the solution of difference scheme
				\eqref{ur03}--\eqref{ur03.1} converges to the solution of
				differential problem \eqref{ur1}--\eqref{ur3} in the grid $L_2$ -
				norm with the rate equal to the order of the approximation error
				$\mathcal{O}(N^{-r_1}+M^{-r_2})$.}
		\end{theorem}		
	
		\section{ A second order numerical differentiation
			formula for the generalized Caputo fractional derivative}
		
		In this section, we construct  a difference analog of the Caputo fractional
		derivative with the  approximation order $\mathcal
		O(\tau^{2})$ and investigate its essential properties.

		Next we consider the uniform grid  $\bar\omega_\tau=\{t_j=j\tau, j=0,
		1, \ldots, M, \tau M=T\}$. Let us find the discrete analog
		of the $\partial_{0t}^{\alpha,\lambda}v(t)$ at the fixed point $t_{j+\sigma}$, $j\in\{0, 1, \ldots, M-1\}$, where $v(t)\in
		\mathcal{C}^3[0,T]$, $\sigma = 1-\alpha/2$. For all $\alpha\in(0,1)$ and $\lambda(t)>0$
		($\lambda'(t)\leq0$, $\lambda(t)\in \mathcal{C}^2[0,T]$) the
		following equalities hold true
		$$
		\partial_{0t_{j+\sigma}}^{\alpha,\lambda(t)}v(t)=
		\frac{1}{\Gamma(1-\alpha)}\int\limits_{0}^{t_{j+\sigma}}\frac{{\lambda(t_{j+\sigma}-\eta)}}{(t_{j+\sigma}-\eta)^\alpha}v'(\eta)d\eta
		$$
		
		$$
		=\frac{1}{\Gamma(1-\alpha)}\left(\sum\limits_{s=1}^{j}\int\limits_{t_{s-1
		}}^{t_{s}}\frac{{\lambda(t_{j+\sigma}-\eta)}}{(t_{j+\sigma}-\eta)^\alpha}v'(\eta)d\eta + \int\limits_{t_j}^{t_{j+\sigma}}\frac{{\lambda(t_{j+\sigma}-\eta)}}{(t_{j+\sigma}-\eta)^\alpha}v'(\eta)d\eta \right)
		$$
		$$
		= \frac{1}{\Gamma(1-\alpha)}\sum\limits_{s=1}^{j}\int\limits_{t_{s-1
		}}^{t_{s}}\frac{{\lambda(t_{j+\sigma}-\eta)}}{(t_{j+\sigma}-\eta)^\alpha}\left(\Pi_{2,s}v(\eta)\right)'d\eta
		$$
		$$
		+\frac{1}{\Gamma(1-\alpha)}
		\sum\limits_{s=1}^{j}\int\limits_{t_{s-1
		}}^{t_{s}}\frac{{\lambda(t_{j+\sigma}-\eta)}}{(t_{j+\sigma}-\eta)^\alpha}\left(v(\eta)-\Pi_{2,s}v(\eta)\right)'d\eta  
		$$
		$$
		+\frac{1}{\Gamma(1-\alpha)}\int\limits_{t_j}^{t_{j+\sigma}}\frac{{\lambda(t_{j+\sigma}-\eta)}}{(t_{j+\sigma}-\eta)^\alpha}\left(\Pi_{1,j}v(\eta)\right)'d\eta 
		$$
		$$+\frac{1}{\Gamma(1-\alpha)} \int\limits_{t_j}^{t_{j+\sigma}}\frac{{\lambda(t_{j+\sigma}-\eta)}}{(t_{j+\sigma}-\eta)^\alpha}\left(v(\eta)-\Pi_{1,j}v(\eta)\right)'d\eta 
		$$
		$$
		=\frac{1}{\Gamma(1-\alpha)}\sum\limits_{s=1}^{j}v_{t, s-1}\int\limits_{t_{s-1
		}}^{t_{s}}\frac{{\lambda(t_{j+\sigma}-\eta)}}{(t_{j+\sigma}-\eta)^\alpha}d\eta
		$$
		$$+
		\frac{1}{\Gamma(1-\alpha)}\sum\limits_{s=1}^{j}v_{\bar tt, s}\int\limits_{t_{s-1
		}}^{t_{s}}\frac{{\lambda(t_{j+\sigma}-\eta)}(\eta-t_{s-1/2})}{(t_{j+\sigma}-\eta)^\alpha}d\eta
		+
		$$
		$$
		+\frac{v_{t,j}}{\Gamma(1-\alpha)}\int\limits_{t_j}^{t_{j+\sigma}}\frac{{\lambda(t_{j+\sigma}-\eta)}}{(t_{j+\sigma}-\eta)^\alpha}d\eta 
		+ R_{1j}^{(1)}+ R_{j j+\sigma}^{(1)} 
		$$

		$$
		=\frac{1}{\Gamma(1-\alpha)}\sum\limits_{s=1}^{j}\left(v_{t, s-1}\int\limits_{t_{s-1
		}}^{t_{s}}\frac{{\lambda_{j-s+\sigma+1/2} -\lambda_{t,j-s+\sigma}(\eta-t_{s-1/2})}}{(t_{j+\sigma}-\eta)^\alpha}d\eta \right.
		$$
		$$
		\left.+\lambda_{j-s+\sigma} v_{\bar tt, s}\int\limits_{t_{s-1
		}}^{t_{s}}\frac{(\eta-t_{s-1/2})}{(t_{j+\sigma}-\eta)^\alpha}d\eta\right)
		$$
		$$
		+\frac{\lambda_{\sigma-1/2} v_{t,j}}{\Gamma(1-\alpha)}\int\limits_{t_j}^{t_{j+\sigma}}\frac{d\eta}{(t_{j+\sigma}-\eta)^\alpha} 
		+ R_{1j}^{(1)} + R_{j j+\sigma}^{(1)}+ R_{1j}^{(2)} + R_{j j+\sigma}^{(2)} + R_{1j}^{(3)}
		$$
		$$
		= \frac{\tau^{1-\alpha}}{\Gamma(2-\alpha)}\sum\limits_{s=1}^{j}\left(v_{t,s-1}(\lambda_{j-s+\sigma+1/2}a_{j-s+1}^{(\alpha)}+(\lambda_{j-s+\sigma}-\lambda_{j-s+\sigma+1})b_{j-s+1}^{(\alpha)})\right.
		$$
		$$
		\left.+
		\lambda_{j-s+\sigma}b_{j-s+1}^{(\alpha)}(v_{t,s}-v_{t,s-1})\right) + \frac{\tau^{1-\alpha}}{\Gamma(2-\alpha)}\lambda_{\sigma-1/2}a_0^{(\alpha)}v_{t,j} + R_{1}^{j+\sigma}
		$$
		$$
		= \frac{\tau^{1-\alpha}}{\Gamma(2-\alpha)}\left((\lambda_{j + \sigma - 1/2} a_{j}^{(\alpha)} - \lambda_{j + \sigma} b_{j}^{(\alpha)})v_{t,0} \right.
		$$
		$$
		+\sum\limits_{s=1}^{j-1}\left(\lambda_{j-s+\sigma-1/2}a_{j-s}^{(\alpha)}+\lambda_{j-s+\sigma}b_{j-s+1}^{(\alpha)}-\lambda_{j-s+\sigma}b_{j-s}^{(\alpha)}\right)v_{t,s}
		$$
		$$
		+ \left.(\lambda_{\sigma-1/2} a_0^{(\alpha)} + \lambda_{\sigma} b_{1}^{(\alpha)})v_{t,j}\right)
		$$
		$$
		=\frac{\tau^{1-\alpha}}{\Gamma(2-\alpha)}\sum\limits_{s=0}^{j}c_{j-s}^{(\alpha)}v_{t,s} + R_{1}^{j+\sigma}.
		$$
		
		where
		$$
		a_0^{(\alpha)} = \sigma^{1-\alpha}, \quad  a_l^{(\alpha)} = (l+\sigma)^{1-\alpha} - (l-1+\sigma)^{1-\alpha},\quad
		$$ 
		$$
		b_l^{(\alpha)}=\frac{1}{2-\alpha}[(l+\sigma)^{2-\alpha}-(l-1+\sigma)^{2-\alpha}]-\frac{1}{2}[(l+\sigma)^{1-\alpha}+(l-1+\sigma)^{1-\alpha}],
		\quad l\geq1,
		$$
		
		$$\lambda_s=\lambda(t_s),\quad
		v_{t,s}=\frac{v(t_{s+1})-v(t_s)}\tau, \quad
		v_{\bar tt,s}=\frac{v(t_{s})-v(t_{s-1})}\tau,
		$$
		$$\Pi_{1,s}v(t)=v(t_{s+1})\frac{t-t_s}{\tau}+v(t_{s})\frac{t_{s+1}-t}{\tau},
		$$
		$$\Pi_{2,s}v(t)=v(t_{s+1})\frac{(t-t_{s-1})(t-t_s)}{2\tau^2}
		$$
		$$
		-
		v(t_{s})\frac{(t-t_{s-1})(t-t_{s+1})}{\tau^2} + v(t_{s-1})\frac{(t-t_{s})(t-t_{s+1})}{2\tau^2},
		$$
		$$
		R_1^{j+\sigma}= R_{1j}^{(1)} + R_{j j+\sigma}^{(1)}+ R_{1j}^{(2)} + R_{j j+\sigma}^{(2)} + R_{1j}^{(3)},
		$$
		$$
		R_{1j}^{(1)} = \frac{1}{\Gamma(1-\alpha)}\sum\limits_{s=1}^{j}\int\limits_{t_{s-1
		}}^{t_{s}}\frac{{\lambda(t_{j+\sigma}-\eta)}}{(t_{j+\sigma}-\eta)^\alpha}\left(v(\eta)-\Pi_{2,s}v(\eta)\right)'d\eta,
		$$
		$$
		R_{j j+\sigma}^{(1)} = \frac{1}{\Gamma(1-\alpha)}\int\limits_{t_j}^{t_{j+\sigma}}\frac{{\lambda(t_{j+\sigma}-\eta)}}{(t_{j+\sigma}-\eta)^\alpha}\left(v(\eta)-\Pi_{1,j}v(\eta)\right)'d\eta,
		$$
		$$
		R_{1j}^{(2)}=\frac{1}{\Gamma(1-\alpha)}\sum\limits_{s=1}^{j}v_{t, s-1}\int\limits_{t_{s-1
		}}^{t_{s}}\frac{{\lambda(t_{j+\sigma}-\eta) - \lambda_{j-s+\sigma+1/2}  +\lambda_{t,j-s+\sigma}(\eta-t_{s-1/2})}}{(t_{j+\sigma}-\eta)^\alpha}d\eta,
		$$
		$$
		R_{j j+\sigma}^{(2)} = \frac{v_{t,j}}{\Gamma(1-\alpha)}\int\limits_{t_j}^{t_{j+\sigma}}\frac{{\lambda(t_{j+\sigma}-\eta) - \lambda_{\sigma-1/2}}}{(t_{j+\sigma}-\eta)^\alpha}d\eta,
		$$
		$$
		R_{1j}^{(3)} = \frac{1}{\Gamma(1-\alpha)}\sum\limits_{s=1}^{j}
		v_{\bar tt, s}\int\limits_{t_{s-1
		}}^{t_{s}}\frac{{(\lambda(t_{j+\sigma}-\eta)-\lambda_{j-s+\sigma})}(\eta-t_{s-1/2})}{(t_{j+\sigma}-\eta)^\alpha}d\eta.
		$$

		Let us consider the below fractional numerical differentiation formula
		for the generalized Caputo fractional derivative of order $\alpha$
		with weighting function $\lambda(t)$ ($0<\alpha<1, \lambda(t)>0,
		\lambda'(t)\leq0$)
		\begin{equation}
			\Delta_{0t_{j+\sigma}}^{\alpha,\lambda(t)}v=\frac{\tau^{1-\alpha}}{\Gamma(2-\alpha)}\sum\limits_{s=0}^{j}c_{j-s}^{(\alpha)}v_{t,s},
			\label{ur4}
		\end{equation}
		where
		$$
		c_0^{(\alpha)} = \lambda_{\sigma-1/2} a_0^{(\alpha)},\quad \text{for}\quad j = 0;\quad \text{and for}\quad j \geq 1, 
		$$
		\begin{equation}
			c_{s}^{(\alpha)}=
			\begin{cases}
				\lambda_{\sigma-1/2} a_0^{(\alpha)} + \lambda_{\sigma} b_{1}^{(\alpha)}, \quad\quad\quad \quad \quad\quad\quad\quad s=0,\\
				\lambda_{s + \sigma-1/2} a_{s}^{(\alpha)} + \lambda_{s + \sigma} b_{s+1}^{(\alpha)} - \lambda_{s + \sigma} b_{s}^{(\alpha)}, \quad\, 1\leq s\leq j-1,\\
				\lambda_{j + \sigma - 1/2} a_{j}^{(\alpha)} - \lambda_{j + \sigma} b_{j}^{(\alpha)},
				\quad\quad\quad\quad\quad\quad\, s=j. \label{FNDF}
			\end{cases}
		\end{equation}
		We call \eqref{ur4} the $_\lambda$L2-1$_\sigma$ - formula for the generalized Caputo fractional derivative.
		
		\begin{lemma}\label{lem_approx} {\it For any $\alpha\in(0,1)$ and $v(t)\in
				\mathcal{C}^3[0,t_{j+1}]$, it is true that
				\begin{equation}
					\partial_{0t_{j+\sigma}}^{\alpha,\lambda(t)}v=\Delta_{0t_{j+\sigma}}^{\alpha,\lambda(t)}v+\mathcal{O}(\tau^{2}), \label{ur5}
				\end{equation}
				where  $\lambda(t)>0$, $\lambda'(t)\leq0$ and $\lambda(t)\in
				\mathcal{C}^{2}[0,t_{j+1}]$.}
		\end{lemma}
		{\bf Proof.} We have
		$\partial_{0t_{j+\sigma}}^{\alpha,\lambda(t)}v-\Delta_{0t_{j+\sigma}}^{\alpha,\lambda(t)}v=R_{1j}^{(1)} + R_{j j+\sigma}^{(1)}+ R_{1j}^{(2)} + R_{j j+\sigma}^{(2)} + R_{1j}^{(3)}$.
		
		Estimate the errors $R_{1j}^{(1)}$,  $R_{j j+\sigma}^{(1)}$, $R_{1j}^{(2)}$, $R_{j j+\sigma}^{(2)}$ and $R_{1j}^{(3)}$:
		$$
		|R_{1j}^{(1)}|=\frac{1}{\Gamma(1-\alpha)}\left|\sum\limits_{s=1}^{j}\int\limits_{t_{s-1
		}}^{t_{s}}\frac{{\lambda(t_{j+\sigma}-\eta)}}{(t_{j+\sigma}-\eta)^\alpha}\left(v(\eta)-\Pi_{2,s}v(\eta)\right)'d\eta\right|
		$$
		$$
		\leq\frac{1}{\Gamma(1-\alpha)}\sum\limits_{s=1}^{j}\left|\,\int\limits_{t_{s-1
		}}^{t_{s}}\left(-\frac{{\lambda'(t_{j+\sigma}-\eta)}}{(t_{j+\sigma}-\eta)^\alpha}+\frac{{\alpha\lambda(t_{j+\sigma}-\eta)}}{(t_{j+\sigma}-\eta)^{\alpha+1}}\right)\left(v(\eta)-\Pi_{2,s}v(\eta)\right)d\eta\right|
		$$
		$$
		\leq\frac{M_3^{j+1}\tau^3}{9\sqrt{3}\Gamma(1-\alpha)}\sum\limits_{s=1}^{j}\,\int\limits_{t_{s-1}}^{t_{s}}\left(\frac{m_1^{j+1}}{(t_{j+\sigma}-\eta)^\alpha}+\frac{{\alpha\lambda(0)}}{(t_{j+\sigma}-\eta)^{\alpha+1}}\right)d\eta
		$$
		$$
		=\frac{M_3^{j+1}\tau^3}{9\sqrt{3}\Gamma(1-\alpha)}\int\limits_{0}^{t_{j}}\left(\frac{m_1^{j+1}}{(t_{j+\sigma}-\eta)^\alpha}+\frac{{\alpha\lambda(0)}}{(t_{j+\sigma}-\eta)^{\alpha+1}}\right)d\eta
		$$
		$$
		\leq \frac{M_3^{j+1}\tau^3}{9\sqrt{3}\Gamma(1-\alpha)}\left(m_1^{j+1}\frac{t_{j+\sigma}^{1-\alpha}}{1-\alpha}+\frac{\lambda(0)}{\sigma^\alpha \tau^\alpha}\right) = \mathcal{O}(\tau^{3-\alpha}),
		$$
		$$
		|R_{j j+\sigma}^{(1)}| = \frac{1}{\Gamma(1-\alpha)}\left|\int\limits_{t_j}^{t_{j+\sigma}}\frac{{\lambda(t_{j+\sigma}-\eta)}}{(t_{j+\sigma}-\eta)^\alpha}\left(v(\eta)-\Pi_{1,j}v(\eta)\right)'d\eta\right|
		$$
		$$
		=\frac{1}{\Gamma(1-\alpha)}\left|\int\limits_{t_j}^{t_{j+\sigma}}\frac{{\lambda(t_{j+\sigma}-\eta)}}{(t_{j+\sigma}-\eta)^\alpha}\left(v'(\eta)-v_{t, j}\right)d\eta\right|
		$$
		$$
		=\left|\frac{v''(t_{j+1})}{\Gamma(1-\alpha)}\int\limits_{t_j}^{t_{j+\sigma}}\frac{{\lambda(t_{j+\sigma}-\eta)}(\eta-t_{j+1/2})}{(t_{j+\sigma}-\eta)^\alpha}d\eta+ \mathcal{O}(\tau^{3-\alpha})\right| 
		$$
		$$
		=\left|\frac{v''(t_{j+1})\lambda(t_{\sigma-1/2})}{\Gamma(1-\alpha)}\int\limits_{t_j}^{t_{j+\sigma}}\frac{{}(\eta-t_{j+1/2})}{(t_{j+\sigma}-\eta)^\alpha}d\eta+ \mathcal{O}(\tau^{3-\alpha})\right| 
		$$
		$$
		=\left|\frac{v''(t_{j+1})\lambda(t_{\sigma-1/2})\sigma^{1-\alpha}\tau^{2-\alpha}}{\Gamma(3-\alpha)}(\sigma - 1 + \alpha/2)+ \mathcal{O}(\tau^{3-\alpha})\right| = \mathcal{O}(\tau^{3-\alpha}),
		$$
		$$
		|R_{1j}^{(2)}| = \frac{1}{\Gamma(1-\alpha)}\left|\sum\limits_{s=1}^{j}v_{t, s-1}\int\limits_{t_{s-1
		}}^{t_{s}}\frac{{\lambda(t_{j+\sigma}-\eta) - \lambda_{j-s+\sigma+1/2}  +\lambda_{t,j-s+\sigma}(\eta-t_{s-1/2})}}{(t_{j+\sigma}-\eta)^\alpha}d\eta\right|
		$$
		$$
		\leq \frac{M_1^{j+1}m_2^{j+1}\tau^2}{4\Gamma(1-\alpha)}\sum\limits_{s=1}^{j}\int\limits_{t_{s-1}}^{t_{s}}\frac{d\eta}{(t_{j+\sigma}-\eta)^\alpha} = 
		\frac{M_1^{j+1}m_2^{j+1}\tau^2}{4\Gamma(1-\alpha)}\int\limits_{0}^{t_{j}}\frac{d\eta}{(t_{j+\sigma}-\eta)^\alpha} 
		$$
		$$
		\leq
		\frac{M_1^{j+1}m_2^{j+1}t_{j+\sigma}^{1-\alpha}\tau^2}{4\Gamma(1-\alpha)}
		= \mathcal{O}(\tau^{2}),
		$$
		$$
		|R_{j j+\sigma}^{(2)}| = \left|\frac{v_{t,j}}{\Gamma(1-\alpha)}\int\limits_{t_j}^{t_{j+\sigma}}\frac{{\lambda(t_{j+\sigma}-\eta) - \lambda_{\sigma-1/2}}}{(t_{j+\sigma}-\eta)^\alpha}d\eta\right|
		$$
		$$
		=\left|\frac{v_{t,j}}{\Gamma(1-\alpha)}\int\limits_{t_j}^{t_{j+\sigma}}\frac{-\lambda'(t_{\sigma-1/2})(\eta-t_{j+1/2}) + \frac{1}{2}\lambda''(\bar \xi)(\eta-t_{j+1/2})^2}{(t_{j+\sigma}-\eta)^\alpha}d\eta\right|
		$$
		$$
		\leq \frac{M_1^{j+1}m_2^{j+1}\sigma^{1-\alpha}}{4\Gamma(2-\alpha)}\tau^{3-\alpha} = \mathcal{O}(\tau^{3-\alpha}),
		$$
		$$
		|R_{1j}^{(3)}| = \frac{1}{\Gamma(1-\alpha)}\left|\sum\limits_{s=1}^{j}
		v_{\bar tt, s}\int\limits_{t_{s-1
		}}^{t_{s}}\frac{{(\lambda(t_{j+\sigma}-\eta)-\lambda_{j-s+\sigma})}(\eta-t_{s-1/2})}{(t_{j+\sigma}-\eta)^\alpha}d\eta\right|
		$$
		$$
		\frac{1}{\Gamma(1-\alpha)}\left|\sum\limits_{s=1}^{j}
		v_{\bar tt, s}\int\limits_{t_{s-1
		}}^{t_{s}}\frac{-\lambda'(\bar \xi_2)(\eta-t_{s})(\eta-t_{s-1/2})}{(t_{j+\sigma}-\eta)^\alpha}d\eta\right|
		$$
		$$
		\leq \frac{M_2^{j+1}m_1^{j+1}\tau^2}{2\Gamma(1-\alpha)}\int\limits_{0
		}^{t_{j}}\frac{d\eta}{(t_{j+\sigma}-\eta)^\alpha}=
		\frac{M_2^{j+1}m_1^{j+1}t_{j+\sigma}^{1-\alpha}\tau^2}{2\Gamma(2-\alpha)} = \mathcal{O}(\tau^{2}) 
		$$
		
		where  $M_k^{j+1}=\max\limits_{0\leq t\leq t_{j+1}}|v^{(k)}(t)|$,
		$m_k^{j+1}=\max\limits_{0\leq t\leq t_{j+1}}|\lambda^{(k)}(t)|$.

		\begin{lemma}\label{lm_pr_1} For all $\alpha\in(0,1)$ and $s=1, 2, 3, \ldots$
			\begin{equation}
				\frac{1-\alpha}{(s+\sigma)^\alpha}<a_s<\frac{1-\alpha}{(s+\sigma-1)^\alpha},
				\label{url2_1}
			\end{equation}
			\begin{equation}
				\frac{\alpha(1-\alpha)}{(s+\sigma+1)^{\alpha+1}}<a_s-a_{s+1}<\frac{\alpha(1-\alpha)}{(s+\sigma-1)^{\alpha+1}},
				\label{url2_2}
			\end{equation}
			\begin{equation}
				\frac{\alpha(1-\alpha)}{12(s+\sigma)^{\alpha+1}}<b_{s}<\frac{\alpha(1-\alpha)}{12(s+\sigma-1)^{\alpha+1}},
				\label{url2_3}
			\end{equation}
		\end{lemma}
		\textbf{Proof.} The validity of Lemma \ref{lm_pr_1} results from the following
			equalities:
			$$
			a_{s}^{(\alpha)}=(1-\alpha)\int\limits_{0}^{1}\frac{d\xi}{(s+\sigma-1+\xi)^{\alpha}},
			$$
			$$
			a_{s}^{(\alpha)}-a_{s+1}^{(\alpha)}=\alpha(1-\alpha)\int\limits_{0}^{1}d\eta\int\limits_{0}^{1}\frac{d\xi}{(s+\sigma-1+\xi+\eta)^{\alpha+1}},
			$$
			$$
			b_{s}^{(\alpha)}=\frac{\alpha(1-\alpha)}{2^{2-\alpha}}\int\limits_{0}^{1}\eta
			d\eta\int\limits_{2(s+\sigma)-1-\eta}^{2(s+\sigma)-1+\eta}\frac{d\xi}{\xi^{\alpha+1}}.
			$$
		
		\begin{lemma}\label{lem_in_JCP} \cite{AlikhanovJCP} For all $\alpha\in(0,1)$ and $s=1, 2, 3, \ldots$
			\begin{equation}
				a_s^{(\alpha)}-b_s^{(\alpha)}>\frac{1-\alpha}{2}(s+\sigma)^{-\alpha},
				\label{lem32_1}
			\end{equation}
			\begin{equation}
				(2\sigma-1)(a_0^{(\alpha)}+b_1^{(\alpha)})-\sigma(a_1^{(\alpha)}+b_2^{(\alpha)}-b_1^{(\alpha)})>\frac{\alpha(1-\alpha)}{4\sigma(1+\sigma)^\alpha}.
				\label{lem32_2}
			\end{equation}
			
		\end{lemma}
		
		\begin{lemma}\label{lem_prop} For any $\alpha\in(0,1)$ and
			$c_s^{(\alpha)}$ ($0\leq s\leq j$, $j\geq 1$) defined in \eqref{FNDF}, the following is valid
			\begin{equation}
				c_j^{(\alpha)} > \frac{1-\alpha}{2}\frac{\lambda_{j + \sigma}}{(j+\sigma)^{\alpha}},
				\label{lem33_1}
			\end{equation}
			\begin{equation}
				(2\sigma-1)c_0^{(\alpha)}-\sigma c_1^{(\alpha)}>0,
				\label{lem33_2}
			\end{equation}
			\begin{equation}
				c_0^{(\alpha)} > c_1^{(\alpha)}>c_2^{(\alpha)}>\ldots >c_{j-1}^{(\alpha)}>c_j^{(\alpha)},
				\label{lem33_3}
			\end{equation}
			where $\sigma = 1 - \alpha/2\in({1}/{2},1)$.
		
		\end{lemma}
		\textbf{Proof.}
			The inequality \eqref{lem33_1} follows from the inequality \eqref{lem32_1} since
			$$
			c_j^{(\alpha)} =  \lambda_{j + \sigma - 1/2} a_{j}^{(\alpha)} - \lambda_{j + \sigma} b_{j}^{(\alpha)}
			$$
			$$
			\geq \lambda_{j + \sigma}(a_{j}^{(\alpha)} - b_{j}^{(\alpha)})>\frac{1-\alpha}{2}\frac{\lambda_{j + \sigma}}{(j+\sigma)^{\alpha}}.
			$$ 
			The inequality \eqref{lem33_2} follows from the inequality \eqref{lem32_2} since
			$$
			(2\sigma-1)c_0^{(\alpha)}-\sigma c_1^{(\alpha)} = (2\sigma-1)(\lambda_{\sigma-1/2} a_0^{(\alpha)} + \lambda_{\sigma} b_{1}^{(\alpha)})
			$$
			$$
			-\sigma (\lambda_{\sigma+1/2} a_{1}^{(\alpha)} + \lambda_{1 + \sigma} b_{2}^{(\alpha)} - \lambda_{1 + \sigma} b_{1}^{(\alpha)})
			$$
			$$
			\geq\lambda_{\sigma}\left((2\sigma-1)(a_0^{(\alpha)}+b_1^{(\alpha)})-\sigma(a_1^{(\alpha)}+b_2^{(\alpha)}-b_1^{(\alpha)})\right)
			$$
			$$
			-\sigma(\lambda_{\sigma}-\lambda_{1+\sigma})b_1^{(\alpha)}
			> \lambda_{\sigma}\frac{\alpha(1-\alpha)}{4\sigma(1+\sigma)^{\alpha}} - (\lambda_{\sigma}-\lambda_{1+\sigma})\frac{\alpha(1-\alpha)}{12\sigma^\alpha}
			$$
			$$
			>(\lambda_{\sigma}-\lambda_{1+\sigma})\frac{\alpha(1-\alpha)}{12\sigma(1+\sigma)^{\alpha}}\left(3-\sigma^{1-\alpha}(1+\sigma)^\alpha\right)>0.
			$$
			The inequality \eqref{lem33_3} for the case $c_0^{(\alpha)}>c_1^{(\alpha)}$ follows from the inequality \eqref {lem33_2}. Let us prove the inequality $c_s^{(\alpha)}>c_{s + 1}^{(\alpha)}$ for $ s = 1, 2, \ldots, j $. The difference $c_s^{(\alpha)} - c_{s + 1}^{(\alpha)} $ satisfies the following estimates
			$$
			c_s^{(\alpha)}-c_{s+1}^{(\alpha)} = \lambda_{s+\sigma-1/2}a_s^{(\alpha)}-\lambda_{s+\sigma+1/2}a_{s+1}^{(\alpha)}-\lambda_{s+\sigma}b_s^{(\alpha)}
			$$
			$$
			+(\lambda_{s+\sigma}+\lambda_{s+\sigma+1})b_{s+1}^{(\alpha)}-\lambda_{s+\sigma+1}b_{s+2}^{(\alpha)}
			$$
			$$
			>\lambda_{s+\sigma}\left(a_s^{(\alpha)}-a_{s+1}^{(\alpha)}-b_s^{(\alpha)}+b_{s+1}^{(\alpha)}\right)+\lambda_{s+\sigma+1}\left(b_{s+1}^{(\alpha)}-b_{s+2}^{(\alpha)}\right)
			$$
			$$
			>\lambda_{s+\sigma}\left(a_s^{(\alpha)}-a_{s+1}^{(\alpha)}-b_s^{(\alpha)}\right)
			$$
			$$
			>\lambda_{s+\sigma}\left(\frac{\alpha(1-\alpha)}{(s+\sigma+1)^{\alpha+1}}-\frac{\alpha(1-\alpha)}{12(s+\sigma-1)^{\alpha+1}}\right)
			$$
			$$
			=\frac{\alpha(1-\alpha)\lambda_{s+\sigma}}{12(s+\sigma+1)^{\alpha+1}}\left(12-\frac{(s+\sigma+1)^{\alpha+1}}{(s+\sigma-1)^{\alpha+1}}\right)>0.
			$$

		\begin{corollary}\label{lem_ineq}  For any function $v(t)$ defined on the grid
			$\overline \omega_{\tau}$ we have the inequality
			\begin{equation}\label{ur0404}
				\left(\sigma v^{j+1} + (1-\sigma)v^j\right)\Delta_{0t_{j+\sigma}}^{\alpha,\lambda(t)}v\geq
				\frac{1}{2}\Delta_{0t_{j+\sigma}}^{\alpha,\lambda(t)}v^2.
			\end{equation}
		\end{corollary}

		\section{ A second order difference scheme for the generalized time-fractional diffusion equation}
		
		Suppose that a solution $u(x,t)\in \mathcal{C}_{x,t}^{4,3}$ of
		problem \eqref{ur1}--\eqref{ur3} exists, and the coefficients of equation
		\eqref{ur1} and the functions $f(x,t)$ and $u_0(x)$ fulfill the
		conditions, necessary for the construction of difference schemes with
		the order of approximation $\mathcal{O}(h^2+\tau^2)$.
		
		Consider the following difference scheme
		\begin{equation}\label{ur6}
			\Delta_{0t_{j+\sigma}}^{\alpha,\lambda(t)}y_i=\Lambda y^{(\sigma)}_i
			+\varphi_i^{j+\sigma}, \quad i=1,2,\ldots,N-1,\quad j=0,1,\ldots,M-1,
		\end{equation}
		\begin{equation}
			y(0,t)=0,\quad y(l,t)=0,\quad t\in \overline \omega_{\tau}, \quad
			y(x,0)=u_0(x),\quad  x\in \overline \omega_{h},\label{ur7}
		\end{equation}
		where
		$$
		\Lambda y_i=\left((ay_{\bar
			x})_x-dy\right)_i
		$$
		$$
		=\frac{a_{i+1}y_{i+1}-(a_{i+1}+a_i)y_i+a_iy_{i-1}}{h^2}-d_iy_i,\quad
		i=1,\ldots,N-1,
		$$
		$y^{j+\sigma} = \sigma y^{j+1} + (1-\sigma)y^j$, 
		$y_{\bar x,i}=(y_i-y_{i-1})/h$,\, $y_{x,i}=(y_{i+1}-y_{i})/h$,
		$a_i^{j+\sigma}=k(x_{i-1/2},t_{j+\sigma})$,\, $d_i^{j+\sigma}=q(x_{i},t_{j+\sigma})$,
		$\varphi_i^{j+\sigma}=f(x_i,t_{j+\sigma})$.
		
		If the solution of problem \eqref{ur1}-\eqref{ur3} $u\in
		\mathcal{C}_{x,t}^{4,3}$ then according to \cite{Samar:77} and the formula 
		\eqref{ur4}, the order of approximation of difference scheme
		\eqref{ur6}--\eqref{ur7} is $\mathcal{O}(h^2+\tau^{2})$.

		\begin{theorem}\label{theor_apr_1} { The difference scheme \eqref{ur6}--\eqref{ur7}
				is unconditionally stable and for its solution the following a
				priory estimate is valid:
				\begin{equation}
					\|y^{j+1}\|_0^2\leq\|y^0\|_0^2+\frac{T^\alpha
						\Gamma(1-\alpha)}{2\lambda(T)c_1}\max\limits_{0\leq j\leq
						M}\|\varphi^j\|_0^2. \label{ur8}
			\end{equation}}
		\end{theorem}
		
		\textbf{Proof.} For the difference operator $\Lambda$ by means of Green's
			first difference formula and the embedding theorem \cite{Samar:77}
			for the functions vanishing at $x=0$ and $x=1$, we arrive at $(-\Lambda
			y,y)\geq 4c_1\|y\|_0^2$, that is for this operator we can
			take $\varkappa=4c_1$.
			
			Considering that difference scheme \eqref{ur6}--\eqref{ur7} has the form
			\eqref{ur03}--\eqref{ur03.1} where $g_s^{j+1}=\frac{\tau^{-\alpha}}{\Gamma(2-\alpha)}c_{j-s}^{(\alpha)}$, then
			Lemma 5 implies validity of the following inequalities:
			$$
			g_0^{j+1}>\frac{{\lambda(t_{j+\sigma})}}{2\Gamma(1-\alpha)t_{j+\sigma}^\alpha}>\frac{{\lambda(T)}}{2\Gamma(1-\alpha)T^\alpha},
			$$
			$$
			\quad g_j^{j+1}>g_{j-1}^{j+1}>\ldots>g_{0}^{j+1},\quad \sigma=1-\alpha/2.
			$$
			Therefore, validity of Theorem \ref{theor_apr_1} follows from Theorem \ref{theor_JCP_1}. 
			
		
		From Theorem \ref{theor_JCP_2} it results that if the solution of problem
		\eqref{ur1}--\eqref{ur3} is sufficiently smooth, the solution of
		difference scheme \eqref{ur6}--\eqref{ur7} converges to the solution
		of the differential problem with the rate equal to the order of the
		approximation error $\mathcal{O}(h^2+\tau^{2})$.
		
		\subsection{Numerical results} 
		
		Numerical computations are carried out for a
		test problem when the function
		$$
		u(x,t)=\sin(\pi x)\left(1+\frac{6-(6+6b t+3b^2t^2+b^3t^3)e^{-b
				t}}{b^4}\right)
		$$
		is the exact solution of problem \eqref{ur1}--\eqref{ur3} with
		$\lambda(t)=e^{-bt}$, $b\geq0$ and the coefficients
		$k(x,t)=2-\cos(xt)$, $q(x,t)=1-\sin(xt)$, $T=1$.

		The errors ($z=y-u$) and convergence order (CO) in the norms
		$\|\cdot\|_0$ and $\|\cdot\|_{\mathcal{C}(\bar\omega_{h\tau})}$,
		where
		$\|y\|_{\mathcal{C}(\bar\omega_{h\tau})}=\max\limits_{(x_i,t_j)\in\bar\omega_{h\tau}}|y|$,
		are shown in Table \ref{tab:table1}.
		
		Table \ref{tab:table1} demonstrates that as the number of the spatial subintervals and
		time steps increases, keeping $3h=\tau$, the maximum error decreases, as it is expected and the convergence
		order of the approximate scheme is
		$\mathcal{O}(h^2)=\mathcal{O}(\tau^{2})$, where the
		convergence order is presented by the formula:
		CO$=\log_{\frac{\tau_1}{\tau_2}}{\frac{\|z_1\|}{\|z_2\|}}$ ($z_{i}$ is the
		error corresponding to $\tau_{i}$).
		
		Table \ref{tab:table2} demonstrates that if $h=1/10000$, then while the number of time steps
		of our approximate scheme is increasing, the maximum error is decreasing, 
		as one can expect and the convergence order of time is
		$\mathcal{O}(\tau^2)$, where the convergence order is presented by the
		following formula:
		CO$=\log_{\frac{\tau_1}{\tau_2}}{\frac{\|z_1\|}{\|z_2\|}}$.

		\begin{table}[h]
			
			\caption{$L_2$ - norm and maximum norm error behavior versus grid size reduction when   $\tau=3h$.}
			\label{tab:table1}
			\begin{tabular}{ccccccc}
				\hline
				$b$ & $\alpha$ & $h$ & {$\max\limits_{0\leq n\leq M}\|z^n\|_0$} & {CO in $\|\cdot\|_0$} & {$\|z\|_{C(\bar \omega_{h \tau})}$} & {CO in $||\cdot||_{C(\bar \omega_{h \tau})}$} \\
				\hline
				1.0 &  0.9 & 1/10  & $4.853172e-4$ &        & $6.860735e-4$ &        \\
				&      & 1/20  & $1.195117e-4$ & 2.0218 & $1.689468e-4$ & 2.0218 \\
				&	   & 1/40  & $2.966661e-5$ & 2.0102 & $4.193765e-5$ & 2.0103 \\
				&	   & 1/80  & $7.407823e-6$ & 2.0017 & $1.047192e-5$ & 2.0017 \\
				&	   & 1/160 & $1.853344e-6$ & 1.9989 & $2.619972e-6$ & 1.9989 \\
				&	   & 1/320 & $4.639354e-7$ & 1.9981 & $6.558408e-7$ & 1.9981 \\
				&	   & 1/640 & $1.161322e-7$ & 1.9982 & $1.641709e-7$ & 1.9981 \\ 	\vspace{2mm}
				&	   & 1/1280& $2.904554e-8$ & 1.9994 & $4.106038e-8$ & 1.9994 \\
				2.0 &  0.5 & 1/10  & $5.695428e-4$ &        & $8.053893e-4$ &         \\
				&      & 1/20  & $1.281254e-4$ & 2.1522 & $1.811924e-4$ & 2.1522  \\
				&	   & 1/40  & $3.111526e-5$ & 2.0419 & $4.387037e-5$ & 2.0462  \\
				&	   & 1/80  & $7.832071e-6$ & 1.9902 & $1.104282e-5$ & 1.9901 \\
				&	   & 1/160 & $1.970207e-6$ & 1.9910 & $2.777898e-6$ & 1.9910 \\
				&	   & 1/320 & $4.952711e-7$ & 1.9921 & $6.983096e-7$ & 1.9921 \\
				&	   & 1/640 & $1.243664e-7$ & 1.9936 & $1.753507e-7$ & 1.9936  \\ \vspace{2mm}
				&	   & 1/1280& $3.125438e-8$ & 1.9925 & $4.406686e-8$ & 1.9925 \\
				
				3.0 &  0.1 & 1/10  & $5.590468e-4$ &        & $7.905373e-4$ &         \\
				&      & 1/20  & $1.378485e-4$ & 2.0199 & $1.949425e-4$ & 2.0198  \\
				&	   & 1/40  & $3.418923e-5$ & 2.0115 & $4.820603e-5$ & 2.0158  \\
				&	   & 1/80  & $8.555678e-6$ & 1.9986 & $1.206419e-5$ & 1.9985 \\
				&	   & 1/160 & $2.140670e-6$ & 1.9988 & $3.018517e-6$ & 1.9988 \\
				&	   & 1/320 & $5.355715e-7$ & 1.9989 & $7.551986e-7$ & 1.9989 \\
				&	   & 1/640 & $1.340154e-7$ & 1.9987 & $1.889726e-7$ & 1.9987  \\
				&	   & 1/1280& $3.349770e-8$ & 2.0003 & $4.723392e-8$ & 2.0003 \\
				\hline
			\end{tabular}
			
		\end{table}	
				
		\begin{table}[h]
			
			\caption{$L_2$ - norm and maximum norm error behavior versus $\tau$-grid size reduction when $h=1/2000$.}
			\label{tab:table2}
			\begin{tabular}{ccccccc}
				\hline
				$b$ & $\alpha$ & $h$ & {$\max\limits_{0\leq n\leq M}\|z^n\|_0$} & {CO in $\|\cdot\|_0$} & {$\|z\|_{C(\bar \omega_{h \tau})}$} & {CO in $||\cdot||_{C(\bar \omega_{h \tau})}$} \\
				\hline
				3.0 &  0.9 & 1/10  & $6.977406e-5$ &        & $9.866179e-5$ &         \\
				&      & 1/20  & $1.700981e-5$ & 2.0363 & $2.405134e-5$ & 2.0364 \\
				&	   & 1/40  & $4.110301e-6$ & 2.0491 & $5.812025e-6$ & 2.0490 \\ \vspace{2mm}
				&	   & 1/80  & $9.171973e-7$ & 2.1639 & $1.297116e-6$ & 2.1637 \\
				
				2.0 &  0.5 & 1/10  & $1.144134e-4$ &        & $1.617383e-4$ &         \\
				&      & 1/20  & $2.825404e-5$ & 2.0177 & $3.994110e-5$ & 2.0177  \\
				&	   & 1/40  & $6.909733e-6$ & 2.0318 & $9.768017e-6$ & 2.0317  \\  \vspace{2mm}
				&	   & 1/80  & $1.621574e-6$ & 2.0912 & $2.292670e-6$ & 2.0910 \\
				
				1.0 &  0.1 & 1/10  & $9.999960e-5$ &        & $1.412912e-4$ &         \\
				&      & 1/20  & $2.495408e-5$ & 2.0026 & $3.525761e-5$ & 2.0027  \\
				&	   & 1/40  & $6.147966e-6$ & 2.0211 & $8.686914e-6$ & 2.0210  \\
				&	   & 1/80  & $1.438581e-6$ & 2.0955 & $2.033257e-6$ & 2.0951 \\
				\hline
			\end{tabular}
		\end{table}

		\newpage
		\section{ A compact difference scheme for the tempered time-fractional diffusion equation.}
		
		In the current section for problem \eqref{ur1}--\eqref{ur3} with a smooth
		solution, we build up a compact difference scheme with the
		approximation order $\mathcal{O}(h^4+\tau^{2})$ for the case
		when  $k=k(t)$ and $q=q(t)$ \cite{Sun3,Sun4}. Next we prove the stability and
		convergence of the constructed difference scheme
		in the grid  $L_2$ - norm with the rate equal to the order of the
		approximation error. The achieved results are
		supported by the numerical computations performed for a test
		example.
		
		To differential problem \eqref{ur1}--\eqref{ur3}, we put into correspondence a difference
		scheme  in the case
		when  $k=k(t)$ and $q=q(t)$:
		\begin{equation}\label{ur10}
			\Delta_{0t_{j+\sigma}}^{\alpha,\lambda(t)}\mathcal{H}_hy_i=a^{j+\sigma}y_{\bar
				xx,i}^{(\sigma)}
			-d^{j+\sigma}\mathcal{H}_hy_i^{(\sigma)}+\mathcal{H}_h\varphi_i^{j+\sigma}, 
		\end{equation}
		\begin{equation}
			y(0,t)=0,\quad y(l,t)=0,\quad t\in \overline \omega_{\tau}, \quad
			y(x,0)=u_0(x),\quad  x\in \overline \omega_{h},\label{ur11}
		\end{equation}
		where $\mathcal{H}_hv_i=v_i+h^2v_{\bar xx,i}/12$, $i=1,\ldots,N-1$,
		$a^{j+\sigma}=k(t_{j+\sigma})$, $d^{j+\sigma}=q(t_{j+\sigma})$,
		$\varphi_i^{j+\sigma}=f(x_i,t_{j+\sigma})$, $y^{j+\sigma}=\sigma y^{j+1} + (1-\sigma)y^j$.
		
		From  \cite{Sun4} and Lemma 2  we deduce that if $u\in
		\mathcal{C}_{x,t}^{6,3}$, then the difference scheme has the
		approximation order  $\mathcal{O}(\tau^{2}+h^4)$.

		\begin{theorem}\label{theor_comp_1}
			The difference scheme
			\eqref{ur10}--\eqref{ur11} is unconditionally stable and for its
			solution the following a priori estimate is valid:
			\begin{equation}\label{ur12}
				\|\mathcal{H}_hy^{j+1}\|_0^2\leq\|\mathcal{H}_hy^0\|_0^2+\frac{T^\alpha \Gamma(1-\alpha)}{\lambda(T)c_1}\max\limits_{0\leq j\leq
					M}\|\mathcal{H}_h\varphi^{j}\|_0^2,
			\end{equation}
		\end{theorem}
		
		\textbf{Proof.} Taking the scalar product of the equation
			\eqref{ur10} with $\mathcal{H}_hy^{(\sigma)}=(\mathcal{H}_hy)^{(\sigma)}$, we
			get
			$$
			(\mathcal{H}_hy^{(\sigma)},\Delta_{0t_{j+\sigma}}^{\alpha,\lambda(t)}\mathcal{H}_hy)-a^{j+\sigma}(\mathcal{H}_hy^{(\sigma)},y_{\bar
				xx}^{(\sigma)})
			$$
			\begin{equation}\label{ur13}
				+d^{j+\sigma}(\mathcal{H}_hy^{(\sigma)},\mathcal{H}_hy^{(\sigma)})=(\mathcal{H}_hy^{(\sigma)},\mathcal{H}_h\varphi^{j+\sigma}).
			\end{equation}
			
			Let us transform the terms in identity \eqref{ur13} as
			$$
			(\mathcal{H}_hy^{(\sigma)},\Delta_{0t_{j+\sigma}}^{\alpha,\lambda(t)}\mathcal{H}_hy)\geq\frac{1}{2}\Delta_{0t_{j+\sigma}}^{\alpha,\lambda(t)}\|\mathcal{H}_hy\|_0^2,
			$$
			$$
			-(\mathcal{H}_hy^{(\sigma)},y_{\bar xx}^{(\sigma)})=-(y^{(\sigma)},y_{\bar
				xx}^{(\sigma)})-\frac{h^2}{12}\|y_{\bar xx}^{(\sigma)}\|_0^2
			$$
			$$
			=\|y_{\bar
				x}^{(\sigma)}]|_0^2-\frac{1}{12}\sum\limits_{i=1}^{N-1}(y_{\bar
				x,i+1}^{(\sigma)}-y_{\bar x,i}^{(\sigma)})^2h
			$$
			$$
			\geq\|y_{\bar x}^{(\sigma)}]|_0^2-\frac{1}{3}\|y_{\bar
				x}^{(\sigma)}]|_0^2=\frac{2}{3}\|y_{\bar
				x}^{(\sigma)}]|_0^2\geq\frac{8}{3}\|y^{(\sigma)}\|_0^2,\quad \text{where}\quad
			\|y]|_0^2=\sum\limits_{i=1}^{N}y_i^2h,
			$$
			$$
			(\mathcal{H}_hy^{(\sigma)},\mathcal{H}_h\varphi^{j+\sigma})\leq\varepsilon\|\mathcal{H}_hy^{(\sigma)}\|_0^2+
			\frac{1}{4\varepsilon}\|\mathcal{H}_h\varphi^{j+\sigma}\|_0^2
			$$
			$$
			=\varepsilon\sum\limits_{i=1}^{N-1}\left(\frac{y_{i-1}^{(\sigma)}+10y_{i}^{(\sigma)}+y_{i+1}^{(\sigma)}}{12}\right)^2h+
			\frac{1}{4\varepsilon}\|\mathcal{H}_h\varphi^{j+\sigma}\|_0^2
			$$
			$$
			\leq\varepsilon\|y^{(\sigma)}\|_0^2+\frac{1}{4\varepsilon}\|\mathcal{H}_h\varphi^{j+\sigma}\|_0^2.
			$$
			Taking into consideration the transformations above, from
			identity \eqref{ur13} with $\varepsilon=\frac{8c_1}{3}$ we get
			the inequality
			$$
			\Delta_{0t_{j+1}}^{\alpha,\lambda(t)}\|\mathcal{H}_hy\|_0^2\leq\frac{1}{8c_1}\|\mathcal{H}_h\varphi^{j+1}\|_0^2.
			$$
			The following procedure is similar to the proof of Theorem 1 in
			\cite{AlikhanovJCP}, and it is left out.
		
		The norm $\|\mathcal{H}_hy\|_0$ is equivalent to the norm $\|y\|_0$,
		which results from the inequalities
		$$
		\frac{5}{12}\|y\|_0^2\leq\|\mathcal{H}_hy\|_0^2\leq\|y\|_0^2.
		$$
		
		Likewise Theorem \ref{theor_JCP_2}, we get the convergence result.
		
		\begin{theorem} { Suppose that
				$u(x,t)\in\mathcal{C}_{x,t}^{6,3}$ is the solution of problem
				\eqref{ur1}--\eqref{ur3} for the case when $k=k(t)$, $q=q(t)$, and
				$\{y_i^j \,|\, 0\leq i\leq N, \, 1\leq j\leq M\}$ is the solution of
				difference scheme \eqref{ur10}--\eqref{ur11}. Then the following holds true
				$$
				\|u(\cdot,t_j)-y^j\|_0\leq C_R\left(\tau^{2}+h^4\right),\quad
				1\leq j\leq M,
				$$
				where $C_R$ is a positive constant not depending on $\tau$ and $h$.}
		\end{theorem}
		
		\subsection{Numerical results} 
		
		In this subsection we present a test
		example for a numerical research of difference scheme
		\eqref{ur10}--\eqref{ur11}.
		
		Examine the following problem:
		\begin{equation}\label{ur14}
			\partial_{0t}^{\alpha,\lambda(t)}u=k(t)\frac{\partial^2u}{\partial x^2}-q(t)u+f(x,t),\,\, 0<x<1,\,\, 0<t\leq 1,
		\end{equation}
		\begin{equation}
			u(0,t)=0,\, u(1,t)=0,\, 0\leq t\leq 1, \, u(x,0)=\sin(\pi x),\,
			0\leq x\leq 1,\label{ur15}
		\end{equation}
		where $\lambda(t)=e^{-bt}$, $b\geq 0$, $ k(t)=2-\sin{(3t)}$, \quad
		$q(t)=1-\cos{(2t)},$
		$$
		f(x,t)=\left[\pi^2g(t)k(t)+g(t)q(t)+\frac{2t^{3-\alpha}e^{-b
				t}}{\Gamma(4-\alpha)}\right]\sin(\pi x),
		$$
		whose exact analytical solution is $u(x,t)=g(t)\sin(\pi x)$,
		where
		$$
		g(t)=1+\frac{2-(2+2b t+b^2t^2)e^{-b t}}{b^3}.
		$$
		
		Table \ref{tab:table3} presents the $L_2$ - norm, the errors of the maximum norm  and the
		time convergence order for $\alpha=0.1, 0.5, 0.9$, where
		$h=1/500$. By this we can see that the time convergence order is
		$2$.
		
		Table \ref{tab:table4} shows the $L_2$ - norm, the maximum norm errors and the
		time convergence order, where
		$\tau=1/2000$. We can see that the order of convergence in space is $4$.
		
		Table \ref{tab:table5} demonstrates that as the number of spatial subintervals and time
		steps increases keeping $\tau=16h^2$, the
		maximum error is reduced, as it is expected, and the convergence order of
		the approximate of the scheme is $\mathcal{O}(h^4+\tau^2)=\mathcal{O}(\tau^2)$.
		\begin{table}[h]
			
			\caption{$L_2$ - norm and maximum norm error behavior compared with $\tau$-grid size reduction  when  $h=1/500$.}
			\label{tab:table3}
			\begin{tabular}{ccccccc}
				\hline
				$b$ & $\alpha$ & $\tau$ & {$\max\limits_{0\leq n\leq M}\|z^n\|_0$} & {CO in $\|\cdot\|_0$} & {$\|z\|_{C(\bar \omega_{h \tau})}$} & {CO in $||\cdot||_{C(\bar \omega_{h \tau})}$} \\
				\hline
				1.0 &  0.9 & 1/10  & $3.870828e-4$ &        & $5.474178e-4$ &         \\
				&      & 1/20  & $9.636762e-5$ & 2.0060 & $1.362844e-4$ & 2.0060 \\
				&	   & 1/40  & $2.398099e-5$ & 2.0066 & $3.391425e-5$ & 2.0066 \\
				&	   & 1/80  & $5.973624e-6$ & 2.0052 & $8.447980e-6$ & 2.0052 \\
				&	   & 1/160 & $1.488446e-6$ & 2.0048 & $2.104980e-6$ & 2.0048 \\ \vspace{2mm}
				&	   & 1/320 & $3.709923e-7$ & 2.0043 & $5.246623e-7$ & 2.0043 \\
				
				2.0 &  0.5 & 1/10  & $1.383725e-4$ &        & $1.956883e-4$ &         \\
				&      & 1/20  & $3.418301e-5$ & 2.0172 & $4.834208e-5$ & 2.0172  \\
				&	   & 1/40  & $8.442745e-6$ & 2.0174 & $1.193984e-5$ & 2.0174  \\
				&	   & 1/80  & $2.092596e-6$ & 2.0124 & $2.959377e-6$ & 2.0124 \\
				&	   & 1/160 & $5.200842e-7$ & 2.0084 & $7.355101e-7$ & 2.0084 \\ \vspace{2mm}
				&	   & 1/320 & $1.295146e-7$ & 2.0056 & $1.831613e-7$ & 2.0056 \\
				
				3.0 &  0.1 & 1/10  & $2.622451e-5$ &        & $3.708705e-5$ &         \\
				&      & 1/20  & $6.094819e-6$ & 2.1052 & $8.619377e-6$ & 2.1052  \\
				&	   & 1/40  & $1.451037e-6$ & 2.0704 & $2.052077e-6$ & 2.0705  \\
				&	   & 1/80  & $3.532982e-7$ & 2.0381 & $4.996392e-7$ & 2.0381 \\
				&	   & 1/160 & $8.699997e-8$ & 2.0217 & $1.230365e-7$ & 2.0218 \\
				&	   & 1/320 & $2.156752e-8$ & 2.0121 & $3.050108e-8$ & 2.0121 \\
				\hline
			\end{tabular}			
		\end{table}	
		
		\begin{table}[h]
			
			\caption{$L_2$ - norm and maximum norm error behavior compared with grid size reduction when $\tau=16h^2$.}
			\label{tab:table4}
			\begin{tabular}{ccccccc}
				\hline
				$b$ & $\alpha$ & $h$ & {$\max\limits_{0\leq n\leq M}\|z^n\|_0$} & {CO in $\|\cdot\|_0$} & {$\|z\|_{C(\bar \omega_{h \tau})}$} & {CO in $||\cdot||_{C(\bar \omega_{h \tau})}$} \\
				\hline
				1.0 &  0.9 & 1/4   & $1.216509e-3$ &        & $1.720403e-3$ &         \\
				&      & 1/8   & $7.463500e-5$ & 4.0267 & $1.055498e-4$ & 4.0267 \\
				&	   & 1/16  & $4.635757e-6$ & 4.0089 & $6.555951e-6$ & 4.0089 \\ \vspace{2mm}
				&	   & 1/32  & $2.818584e-7$ & 4.0397 & $3.986080e-7$ & 4.0397 \\
				
				2.0 &  0.5 & 1/4   & $1.133742e-3$ &        & $1.603353e-3$ &         \\
				&      & 1/8   & $6.956352e-5$ & 4.0266 & $9.837767e-4$ & 4.0266 \\
				&	   & 1/16  & $4.327824e-6$ & 4.0066 & $6.120468e-6$ & 4.0066 \\ \vspace{2mm}
				&	   & 1/32  & $2.702171e-7$ & 4.0014 & $3.821448e-7$ & 4.0014 \\
				
				3.0 &  0.1 & 1/4   & $1.086389e-3$ &        & $1.536387e-3$ &         \\
				&      & 1/8   & $6.666005e-5$ & 4.0265 & $9.427155e-4$ & 4.0266 \\
				&	   & 1/16  & $4.147156e-6$ & 4.0066 & $5.864965e-6$ & 4.0066 \\
				&	   & 1/32  & $2.588975e-7$ & 4.0016 & $3.661364e-7$ & 4.0016 \\
				\hline
			\end{tabular}
			
		\end{table}

		\begin{table}[h]
			
			\caption{$L_2$ - norm and maximum norm error behavior compared with the grid size reduction when $\tau=16h^2$.}
			\label{tab:table5}
			\begin{tabular}{ccccccc}
				\hline
				$b$ & $\alpha$ & $\tau$ & {$\max\limits_{0\leq n\leq M}\|z^n\|_0$} & {CO in $\|\cdot\|_0$} & {$\|z\|_{C(\bar \omega_{h \tau})}$} & {CO in $||\cdot||_{C(\bar \omega_{h \tau})}$} \\
				\hline
				1.0 &  0.9 & 1/10   & $3.828076e-4$ &        & $5.413717e-4$ &         \\
				&      & 1/20   & $9.462480e-5$ & 2.0163 & $1.362844e-4$ & 2.0163 \\
				&	   & 1/40   & $2.352703e-5$ & 2.0078 & $3.327224e-5$ & 2.0079 \\
				&	   & 1/80   & $5.807158e-6$ & 2.0184 & $8.212562e-6$ & 2.0184 \\
				&	   & 1/160  & $1.450182e-6$ & 2.0015 & $2.050867e-6$ & 2.0016 \\
				&	   & 1/320  & $3.605780e-7$ & 2.0078 & $5.099343e-7$ & 2.0078 \\
				&	   & 1/640  & $9.010072e-8$ & 2.0007 & $1.274216e-7$ & 2.0007 \\
				&	   & 1/1280 & $2.241364e-8$ & 2.0071 & $3.169767e-8$ & 2.0072 \\ \vspace{2mm}
				&	   & 1/2560 & $5.591086e-9$ & 2.0031 & $7.906995e-9$ & 2.0032 \\
				
				2.0 &  0.5 & 1/10   & $1.342903e-4$ &        & $1.899152e-4$ &         \\
				&      & 1/20   & $3.253876e-5$ & 2.0451 & $4.601676e-5$ & 2.0451 \\
				&	   & 1/40   & $8.015256e-6$ & 2.0213 & $1.133528e-5$ & 2.0213 \\
				&	   & 1/80   & $1.935905e-6$ & 2.0497 & $2.737783e-6$ & 2.0497 \\
				&	   & 1/160  & $4.839828e-7$ & 2.0000 & $6.844551e-7$ & 2.0000 \\
				&	   & 1/320  & $1.196592e-7$ & 2.0160 & $1.692237e-7$ & 2.0160 \\
				&	   & 1/640  & $3.002070e-8$ & 1.9949 & $4.245569e-8$ & 1.9949 \\
				&	   & 1/1280 & $7.438279e-9$ & 2.0129 & $1.051931e-8$ & 2.0129 \\ \vspace{2mm}
				&	   & 1/2560 & $1.857629e-9$ & 2.0015 & $2.627084e-9$ & 2.0015 \\
				
				3.0 &  0.1 & 1/10   & $2.218725e-5$ &        & $3.137751e-5$ &         \\
				&      & 1/20   & $4.434359e-6$ & 2.3229 & $6.271131e-6$ & 2.3229 \\
				&	   & 1/40   & $1.019302e-6$ & 2.1211 & $1.441511e-6$ & 2.1211 \\
				&	   & 1/80   & $3.005858e-7$ & 1.7617 & $4.250925e-7$ & 1.7617 \\
				&	   & 1/160  & $7.429821e-8$ & 2.0163 & $1.050735e-7$ & 2.0163 \\
				&	   & 1/320  & $1.967234e-8$ & 1.9171 & $2.782089e-8$ & 1.9171 \\
				&	   & 1/640  & $4.756649e-9$ & 2.0481 & $6.726918e-9$ & 2.0481 \\
				&	   & 1/1280 & $1.243872e-9$ & 1.9351 & $1.759102e-9$ & 1.9351 \\
				&	   & 1/2560 & $3.110283e-10$& 1.9997 & $4.398604e-10$& 1.9997 \\
				\hline
			\end{tabular}
		\end{table}	
	
		\section{Conclusion}
		
		In the current paper, we study the stability and convergence of a difference schemes
		which approximate the time fractional diffusion equation with
		generalized memory kernel. We have built a new difference approximation
		of the generalized Caputo fractional derivative with the
		approximation order $\mathcal{O}(\tau^{2})$.
		The essential features of this difference operator are investigated.
		We have also constructed some new difference schemes of the second and fourth approximation order
		in space and the second approximation order in time for the
		generalized time fractional diffusion equation with variable
		coefficients. The stability and convergence
		of these schemes in the grid $L_2$ - norm with the rate equal to the
		order of the approximation error are proven as well. The method can be
		without difficulty expanded to other time fractional partial differential
		equations with any other boundary conditions.
		
		Numerical tests thoroughly confirming the achieved theoretical
		results are implemented. In all the computations Julia v1.6.2 is used.
		\vskip 5mm
		\textbf{Funding. }This research was jointly funded by Russian Foundation for Basic Research (RFBR) and Natural  Science  Foundation of China (NSFC), grant numbers 20-51-53007 and 12011530058. The Russian Foundation for Basic Research (RFBR), grant number 19-31-90094, also supported this work.


\newpage

	\begin{thebibliography}{5}
		
		\bibitem{OldSpan}
		Oldham, K.B.; Spanier, J. \textit{The Fractional Calculus}; Academic Press:
		New York, USA, 1974.
		
		\bibitem{Podl}
		Podlubny, I. \textit{Fractional Differential Equations}; Academic Press: San
		Diego, USA, 1999.
		
		\bibitem{Hilfer}
		Hilfer, R. \textit{Applications of Fractional Calculus in Physics}; World
		Scientific: Singapore, 2000.
		
		\bibitem{KilbSrivTruj}
		Kilbas, A.A.; Srivastava, H.M.; Trujillo, J.J. \textit{Theory and Applications
			of Fractional Differential Equation}; Elsevier: Amsterdam, Netherlands, 2006.
		
		\bibitem{Fokk_Plan_Smol}
		Sandev, T.; Chechkin, A.; Kantz, H.; Metzler, R. Diffusion and
		Fokker-Planck-Smoluchowski equations with generalized memory kernel.
		{\em Fract. Calc. Appl. Anal.} {\bf 2015}, {\em 18}, 1006--1038.
		
		\bibitem{Alikh:10}
		Alikhanov, A.A. A priori estimates for solutions of boundary value problems for fractional-order equations. {\em Differ. Equ.} {\bf 2010}, {\em 46},
		660--666.
		
		\bibitem{Alikh:12}
		A.A. Alikhanov, Boundary value problems for the diffusion equation of the variable order in differential and difference settings. {\em Appl. Math. Comput.} {\bf 2012}, {\em 219}, 3938--3946.
		
		\bibitem{AlikhanovJCP}
		Alikhanov, A. A.  A new difference scheme for the time fractional
		diffusion equation. {\em J. Comput. Phys.} {\bf 2015}, {\em 280}, 424--438.
		
		\bibitem{Alikh_15}
		Alikhanov, A.A. Numerical methods of solutions of boundary value
		problems for the multi-term variable-distributed order diffusion
		equation. {\em Appl. Math. Comput.} {\bf 2015}, {\em 268}, 12--22. 
		
		\bibitem{Alikh_16}
		Alikhanov, A. A. Stability and convergence of difference schemes for
		boundary value problems for the fractional-order diffusion equation. 
		{\em Comput. Math. and Math. Phys.} {\bf 2016}, {\em 56}, 561--575.
		
		\bibitem{Alikh_17_gen}
		Alikhanov, A. A. A time-fractional diffusion equation with generalized memory kernel in differential and difference settings with smooth solutions. {\em Comput. Methods Appl. Math.}  {\bf 2017}, {\em 17}, 647--660.
		
		\bibitem{Alikh_17}
		Gao, G.-H.; Alikhanov, A.A. ; Sun Z.-Z. The Temporal Second Order
		Difference Schemes Based on the Interpolation Approximation for
		Solving the Time Multi-term and Distributed-Order Fractional
		Sub-diffusion Equations. {\em J. Sci. Comput.} {\bf 2017}, {\em 73},  93--121. 
		
		\bibitem{Khibiev}
		Khibiev, A. Kh. Stability and convergence of difference schemes for the multi-term time-fractional diffusion equation with generalized memory kernels. {\em J. Samara State Tech. Univ., Ser. Phys. Math. Sci.}, {\bf 2019}, {\em 23},  582--597.
		
		
		\bibitem{ShkhTau:06}
		Shkhanukov-Lafishev, M. Kh.; Taukenova, F.I. Difference methods for
		solving boundary value problems for fractional differential
		equations. {\em Comput. Math. and Math. Phys.} {\bf 2006}, {\em 46}, 785--1795.
		
		\bibitem{SakaYama}
		Sakamoto, K.; Yamamoto, M. Initial value/boundary value problems for
		fractional diffusion-wave equations and applications to some inverse
		problems. {\em J. Math. Anal. Appl.} {\bf 2011} {\em 382} 426--447.
		
		\bibitem{Luch}
		Luchko, Y. Initial-boundary-value problems for the one-dimensional
		time-fractional diffusion equation. {\em Fract. Calc. Appl. Anal.} {\bf 2012}, {\em 15},  141--160.
		
		\bibitem{Alikh_17_1}
		Alikhanov, A.A. A Difference Method for Solving the Steklov
		Nonlocal Boundary Value Problem of Second Kind for the
		Time-Fractional Diffusion Equation.  {\em Comput. Methods Appl. Math.}  {\bf 2017}, {\em 17}, 1--16.
		
		
		\bibitem{Stynes}
		Stynes, M.; O'Riordan, E.; Gracia, J. L. Error analysis of a finite
		difference method on graded meshes for a time-fractional diffusion
		equation.  {\em SIAM J. Numer. Anal.} {\bf 2016}, {\em 55}, 1057--1079.
		
		
		\bibitem{Lazarov_1}
		Jin, B.; Lazarov, R.; Zhou Z. An analysis of the L1 scheme for the
		subdiffusion equation with nonsmooth data. {\em IMA J. Numer. Anal.} {\bf 2015}, {\em 36},
		197--221.
		
		\bibitem{Lazarov_2}
		Jin, B.; Lazarov, R.; Sheen, D.; Zhou, Z. Error estimates for
		approximations of distributed order time fractional diffusion with
		nonsmooth data. {\em Fract. Calc. Appl. Anal.} {\bf 2015}, {\em 19}, 69--93.
		
		\bibitem{Sun6}
		Gao, G. H.; Sun, Z. Z.; Zhang, H. W. A new fractional numerical
		differentiation formula to approximate the Caputo fractional
		derivative and its applications. {\em J. Comput. Phys.} {\bf 2014}, {\em 259}, 33--50.
		
		\bibitem{Samar:77}
		Samarskii, A. A.; \textit{The Theory of Difference Schemes}; Marcel Dekker Inc.: New York, USA, 2001; p. 762.
		
		\bibitem{Sun3}
		Du, R.; Cao, W. R.; Sun, Z. Z. A compact difference scheme for the
		fractional diffusion-wave equation. {\em Appl. Math. Model.} {\bf 2010}, {\em 34}, 2998--3007.
		
		\bibitem{Sun4}
		Gao, G. H.; Sun, Z. Z. A compact difference scheme for the fractional
		subdiffusion equations. {\em J. Comput. Phys.} {\bf 2011}, {\em 230}, 586--595.
		
	\end{thebibliography}
\end{document}